\documentclass[11pt]{amsart}
\usepackage{graphicx}
\usepackage{graphics}
\usepackage{amsmath}
\usepackage{amscd}
\usepackage{latexsym,url}
\usepackage[all]{xy}
\begin{document}
\textwidth 5.5in
\textheight 8.3in
\evensidemargin .75in
\oddsidemargin.75in

\newtheorem{lem}{Lemma}[section]
\newtheorem{conj}[lem]{Conjecture}
\newtheorem{defn}[lem]{Definition}
\newtheorem{thm}[lem]{Theorem}
\newtheorem{cor}[lem]{Corollary}
\newtheorem{prob}[lem]{Problem}
\newtheorem{exm}[lem]{Example}
\newtheorem{rmk}[lem]{Remark}
\newtheorem{que}[lem]{Question}
\newtheorem{prop}[lem]{Proposition}
\newtheorem{clm}[lem]{Claim}
\newcommand{\p}[3]{\Phi_{p,#1}^{#2}(#3)}
\def\Z{\mathbb Z}
\def\R{\mathbb R}
\def\g{\overline{g}}
\def\odots{\reflectbox{\text{$\ddots$}}}
\newcommand{\tg}{\overline{g}}
\def\proof{{\bf Proof. }}
\def\ee{\epsilon_1'}
\def\ef{\epsilon_2'}
\title{Boundary-sum irreducible finite order corks}
\author{Motoo Tange}
\thanks{The author was partially supported by JSPS KAKENHI Grant Number 17K14180}
\address{Institute of Mathematics, University of Tsukuba,
 1-1-1 Tennodai, Tsukuba, Ibaraki 305-8571, Japan}
\email{tange@math.tsukuba.ac.jp}
\date{\today}
\maketitle
\begin{abstract}
We prove for any positive integer $n$ there exist boundary-sum irreducible ${\mathbb Z}_n$-corks with Stein structure.
Here `boundary-sum irreducible' means the manifold is indecomposable with respect to boundary-sum.
We also verify that some of the finite order corks admit hyperbolic boundary by HIKMOT.
\end{abstract}
%
%
%
%

\section{Introduction}
\label{intro}
\subsection{$G$-corks.}
\begin{figure}[htpb]
\begin{center}
\includegraphics{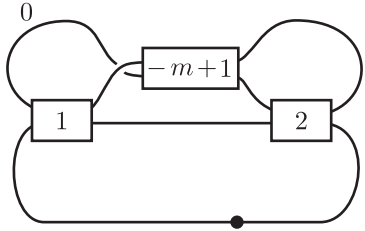}
\caption{The handle diagram of $C(m)$.}
\label{AKB}
\end{center}
\end{figure}

{\it Cork} $(C,g)$ is a pair of a compact contractible (Stein\footnote{This condition is included in some original papers by Akbulut et.al., for example \cite{AK}}) 4-manifold $C$ and a diffeomorphism $g$ on the boundary $\partial C$
that $g$ cannot extend to the inside $C$ as a smooth diffeomorphism.
{\it Cork twist} means the 4-dimensional surgery by the following cut-and-paste
$$X'=(X-C)\cup_gC.$$
The manifold presented by the diagram as in {\sc Figure}~\ref{AKB} becomes a cork.
The map $g$ is the 180$^\circ$ rotation about the horizontal line in the picture.
In particular, $C(1)$ is the first cork which was used by Akbulut.
Here a box with the integer $x$ stands for the $x$-fold right handed full twist.

In the definition of the original cork the condition $g^2=\text{id}_{\partial C}$ is included.
Recently, in some papers the order of gluing map $g$ is generalized to finite order (\cite{TM1}, \cite{AKMR}), infinite order \cite{G} or generally any group $G$ in \cite{AKMR}.
In terms of the view by Auckly, Kim, Melvin, and Ruberman \cite{AKMR},
if a group $G$ smoothly and effectively acts on the boundary of a contractible 4-manifold $C$ and any non-trivial diffeomorphism $g\in G$ cannot smoothly extend to the inside $C$, then the pair $(C,G)$ is called a {\it $G$-cork}.

\begin{figure}[htbp]
\begin{center}
\includegraphics[width=.5\textwidth]{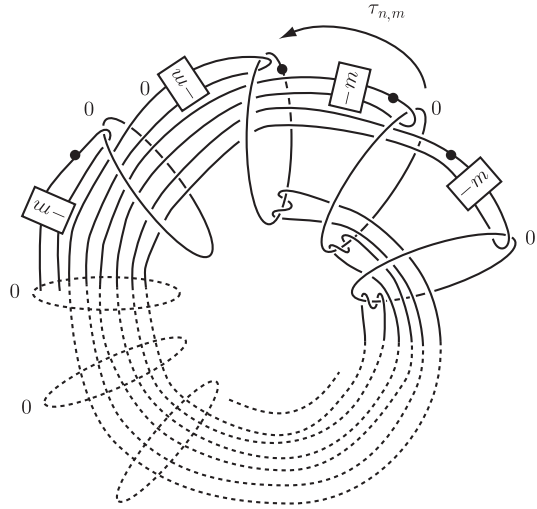}
\caption{The handle decomposition of $C_{n,m}$.}
\label{Cnm}
\end{center}
\end{figure}

\begin{figure}[htbp]
\begin{center}
\includegraphics[width=.5\textwidth]{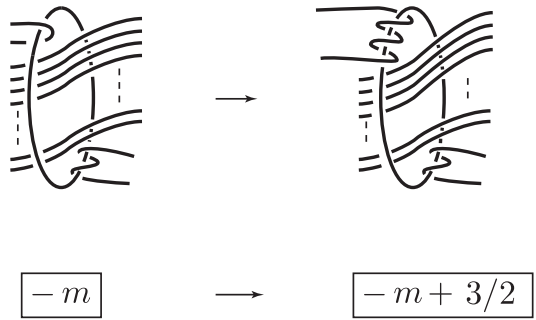}
\caption{The modification.}
\label{beisotopy}
\end{center}
\end{figure}

As examples of finite order corks, the author \cite{TM1} gave pairs of $(X_{n,m},\tau_{n,m})$ for $X=C, D, E$ and or generally, $X=X({\bf x})$ for $\{\ast,0\}$-sequence ${\bf x}\neq (0,\cdots, 0)$ or $(\ast,\cdots, \ast)$, where we call such a sequence ${\bf x}$ {\it non-trivial}.
The diffeomorphism $\tau_{n,m}$ is the $2\pi/n$-rotation with respect to the diagram.
In the paper \cite{TM1}, we put the index $X$ on $\tau_{n,m}$, like the notation $\tau^X_{n,m}$.
We remove the indexes if it is understood from the context.
We describe $C_{n,m}$ in {\sc Figure}~\ref{Cnm}.
$D_{n,m}$ is obtained by exchanging all the dots and $0$-framings in $C_{n,m}$.
$E_{n,m}$ is obtained by modifying $C_{n,m}$ as in {\sc Figure}~\ref{beisotopy}.
The concrete diagrams for these examples are described in \cite{TM1}.

\begin{thm}[\cite{TM1}]
\label{tange}
For $X=C,D, E$ or $X({\bf x})$, for any non-trivial $\{\ast,0\}$-sequence ${\bf x}$, $(X_{n,m},\tau_{n,m})$ is a finite order cork.
Furthermore, $C_{n,m}$ is a ${\mathbb Z}_n$-cork with Stein structure.
\end{thm}

Auckly, Kim, Melvin, and Ruberman \cite{AKMR} gave the examples of $G$-corks for any finite subgroup $G$ of $SO(4)$.
\begin{thm}[\cite{AKMR}]
Let $G$ be any finite subgroup in $SO(4)$.
Then there exists a $G$-cork.
\end{thm}
Let $Y_1,Y_2$ be two $n$-manifolds with boundary.
We call the surgery of attaching an $n$-dimensional 1-handle along two neighborhoods of $p_i\in \partial Y_i$ {\it boundary-sum}
and the resulting manifold as $X_1\natural X_2$.
Their Stein corks in \cite{AKMR} were constructed by the boundary-sum of several copies of $C(1)$.
They also announce the existence of finite order cork with hyperbolic boundary in \cite{AKMR}.

We say that an $n$-manifold $X$ with boundary is {\it boundary-sum irreducible} if $X=X_1\natural X_2$, then $X_1$ or $X_2$ is homeomorphic to an $n$-disk.
If $X$ is not boundary-sum irreducible, then we call $X$ {\it boundary-sum reducible}.
Here a 4-manifold $X$ is called {\it irreducible} if for any connected-sum decomposition $X=X_1\#X_2$, $X_1$ or $X_2$ is a homotopy $n$-sphere.
We call a 3-manifold $Y$ {\it prime} if for any connected-sum decomposition $Y=Y_1\#Y_2$, $Y_1$ or $Y_2$ is a $3$-sphere.
The following holds.
\begin{lem}
\label{equiv}
Let $X$ be a 4-manifold.
If $X$ is irreducible and $\partial X$ is prime, then $X$ is boundary-sum irreducible.
\end{lem}

The problem of whether the examples $X_{n,m}$ in Theorem~\ref{tange} are boundary-sum irreducible corks or not has remained.
Our main theorem answers this question for the case of $X({\bf x})_{n,m}$.
\begin{thm}
\label{main}
For any integer $m$ and positive integer $n$.
There exist boundary-sum irreducible ${\mathbb Z}_n$-corks $(C_{n,m},\tau_{n,m})$ with Stein structure.
For any non-trivial $\{\ast, 0\}$-sequence ${\bf x}$, $(X({\bf x})_{n,m},\tau_{n,m})$ are boundary-sum irreducible 
finite order corks.

Another variation $(E_{n,1},\tau_{n,1})$ is boundary-sum irreducible ${\mathbb Z}_n$-corks.
\end{thm}
Indeed, $X({\bf x})_{n,m}$ is irreducible and $\partial X({\bf x})_{n,m}$ is prime.
This result means that $(X({\bf x})_{n,m},\tau_{n,m})$ is a different finite order cork from the one used in  Theorem A in \cite{AKMR}.
We do not know whether our examples are different from their finite order corks with hyperbolic boundary.

We can show the following result which follows immediately from the proof of Theorem~\ref{main}.
We set $Y_{n,m}:=\partial C_{n,m}$.
Clearly this 3-manifold is diffeomorphic to any $\partial X_{n,m}({\bf x})$, for any $\{\ast,0\}$-sequence ${\bf x}$.
We set $Y'_{n,m}:=\partial E_{n,m}$.
\begin{thm}
\label{irre}
Let $n,m$ be integers as above.
Then $Y_{n,m}$ and $Y'_{n,1}$ are prime homology spheres.
\end{thm}
Furthermore we can prove the following hyperbolicity.
In \cite{AKMR} they suggested any $Y_{n,m}$ would be a hyperbolic 3-manifold.
\begin{thm}
\label{hyper}
Let $n,m$ be integers with $0\le m\le 2$ and $1\le n\le 4$.
$Y_{n,m}$ and $Y'_{n,m}$ are hyperbolic 3-manifolds.
\end{thm}
These are direct results by the computer software HIKMOT \cite{HIKMOT2}.
It is proven that $Y_{1,m}=Y'_{1,m}=\partial C(m)$ are hyperbolic 3-manifolds in \cite{KOU}, by using the fact that these are Dehn surgeries of the pretzel knot $Pr(-3,3,-3)$.
We put a question here.
\begin{que}
Let $X$ be $X({\bf x})$ for non-trivial $\{\ast,0\}$-sequence or $E$.
\begin{itemize}
\item Is $(X_{n,m},\tau_{n,m})$ finite order cork with Stein structure?
\item Is $(X_{n,m},\tau_{n,m})$ finite order cork with hyperbolic boundary?
\end{itemize}
\end{que}
Notice that it is not known at all what reflection for the exotic structures does a cork twist by a cork with hyperbolic boundary give.
This theme is left up to a future study of exotic 4-manifolds.
\section*{Acknowledgements}
I thank for Akitoshi Takayasu and Hidetoshi Masai for helps of the installation of HIKMOT and some useful comments and some strategies to find hyperbolic solution.
I also thank for Kouichi Yasui, Kouki Sato, Takahiro Oba, and Robert Gompf for giving me useful advice, comments and suggestion.

\section{Primeness of $K_{n,m}$ and $K'_{n,m}$.}
\begin{figure}[htbp]
\begin{center}
\includegraphics{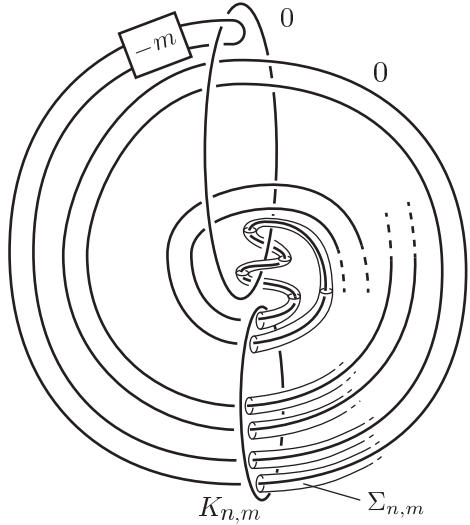}
\caption{$K_{n,m}$ in $Y(m)$}
\label{seifert}
\end{center}
\end{figure}

\begin{figure}[htbp]
\begin{center}
\includegraphics{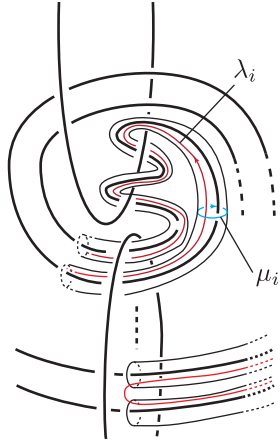}
\caption{Generators of $H_1(\Sigma_{n,m})$.}
\label{generators}
\end{center}
\end{figure}

\begin{figure}[htbp]
\begin{center}
\includegraphics{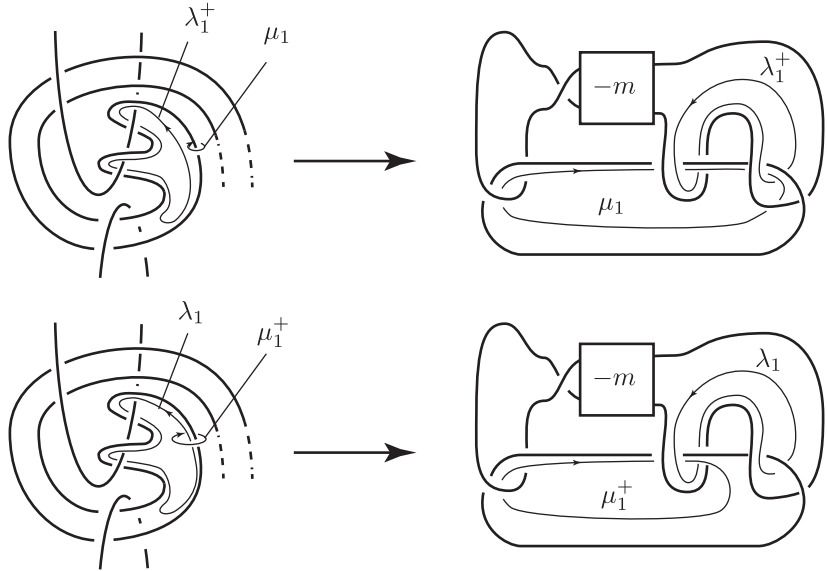}
\caption{$lk(\lambda_i^+,\mu_j)$ and $lk(\lambda_i^+,\mu_j)$}
\label{linking1}
\end{center}
\end{figure}

$Y_{n,m}$ and $Y'_{n,m}$ are $n$-fold cyclic branched covers of $Y(m):=\partial C(m)$ with the branch locus $K_{n,m}$ and $K'_{n,m}$ respectively.
See {\sc Figure}~\ref{seifert} for $K_{n,m}$.
The picture of $K'_{n,m}$ is obtained by modifying the diagram of $K_{n,m}\subset Y(m)$ in {\sc Figure}~\ref{seifert} according to {\sc Figure}~\ref{beisotopy}.
In this picture, the slice disks of $K_{n,m}$ and $K'_{n,m}$ intersect with the 0-framed 2-handle at $2n$ points.
Let $d(K)$ be the top degree of the symmetrized Alexander polynomial $\Delta_K(t)$.
\begin{lem}
For any integer $m$ and positive integer $n$, the Alexander polynomials of $K_{n,m}$ and $K'_{n,m}$ are $\Delta_{K_{n,m}}\doteq 2t^{n}-5+2t^{-n}$ and $\Delta_{K'_{n,m}}=6t^{n}-13+6t^{-n}$. 
Furthermore, the genera of $K_{n,m}$ and $K'_{n,m}$ are $n$.
\end{lem}
Note that the computation in the case of $K_{1,m}$ was done in \cite{KOU}.\\
\begin{proof}
$K_{n,m}$ has the genus $n$ Seifert surface $\Sigma_{n,m}$ as in {\sc Figure}~\ref{seifert}.
We compute the Seifert matrix for $\Sigma_{n,m}$.
We take the generators $\{\lambda_i,\mu_i|i=1,\cdots,n\}$ in $H_1(\Sigma_{n,m})$ as in {\sc Figure}~\ref{generators}.

We define $\lambda_i^+$ and $\mu_i^+$ to be the parallel transforms in the 
one side of the neighborhood of $\Sigma_{n,m}$.
Consider the order of the generators as 
$$\lambda_1,\lambda_2,\cdots,\lambda_n,\mu_1,\mu_2,\cdots,\mu_n.$$
The $(r,s)$-entry of the Seifert matrix $S_{n,m}$ is the linking number $lk(x_r^+,x_s)$,
where $x_i$ is the $i$-th generator above.
Here we have the following:
$$lk(\lambda_i^+,\lambda_j)=0,\ lk(\mu_i^+,\mu_j)=0,\ \ lk(\lambda_i^+,\mu_j)=\begin{cases}2&i\le j\\1&i>j\end{cases}$$
and
$$lk(\mu_i^+,\lambda_j)=\begin{cases}1&i\le j\\2&i>j.\end{cases}$$
These calculations are done by considering the linking indicated in {\sc Figure}~\ref{linking1}. 
The Seifert matrix $S_n$ is $\begin{pmatrix}O_n&A_n\\B_n&O_n\end{pmatrix}$, where $O_n$ is the $n\times n$ zero matrix, $A_n$ and $B_n$ are $n\times n$ matrices satisfying the following:
$$A_n=(a_{ij}),\  a_{ij}=\begin{cases}2&j\ge i\\1&j< i.\end{cases}\text{ and }B_n=(b_{ij}),\ b_{ij}=\begin{cases}1&j\ge i\\2&j< i.\end{cases}$$
Then we have
\begin{eqnarray*}
\Delta_{K_{n,m}}&=&\det(tS_n-S_n^T)=\det\begin{pmatrix}O_n&tA_n-B_n^T\\tB_n-A_n^T&O_n\end{pmatrix}\\
&=&(-1)^n\det(tA_n-B_n^T)\det(tB_n-A_n^T)\\
&=&\det(tA_n-B_n^T)\det(A_n-tB_n^T).
\end{eqnarray*}
We set
$(\alpha_{ij})=tA_n-B_n^T$, where $\alpha_{ij}=\begin{cases}2t-2&j>i\\2t-1&i=j\\t-1&j<i.\end{cases}$
We define $\det(tA_n-B_n^T)$ to be $\alpha_n$.
By expanding $\alpha_n$ and deforming it, we have
$$\alpha_n=\det\begin{pmatrix}1&2t-2&\cdots&\cdots&2t-2\\-t&2t-1&2t-2&\cdots&2t-2\\0&t-1&2t-1&\ddots&\vdots\\\vdots&\vdots&\ddots&\ddots&2t-2\\0&t-1&\cdots&t-1&2t-1\end{pmatrix}=\alpha_{n-1}+t\beta_{n-1},$$
where $\beta_{n-1}$ is the $(n-1)\times (n-1)$ matrix satisfying the following:
\begin{eqnarray*}
\beta_{n-1}&=&
\det
\begin{pmatrix}
2t-2&2t-2&\cdots&\cdots&2t-2\\
t-1&2t-1&2t-2&\cdots&2t-2\\
\vdots&t-1&2t-1&\ddots&\vdots\\
\vdots&\vdots&\ddots&\ddots&2t-2\\
t-1&t-1&\cdots&t-1&2t-1\end{pmatrix}\\
&=&
\det
\begin{pmatrix}
0&2t-2&\cdots&\cdots&2t-2\\
-t&2t-1&2t-2&\cdots&2t-2\\
0&t-1&2t-1&\ddots&\vdots\\
\vdots&\vdots&\ddots&\ddots&2t-2\\
0&t-1&\cdots&t-1&2t-1\end{pmatrix}=t\beta_{n-2}
\end{eqnarray*}
$$\beta_{2}=\det\begin{pmatrix}2t-2&2t-2\\t-1&2t-1\end{pmatrix}=2t(t-1).$$
Thus $\beta_{n-1}=2t^{n-2}(t-1)$, therefore, we have
$\alpha_n=2t^2-1+\sum_{k=3}^n2t^{k-1}(t-1)=2t^n-1$.

By using the following equality
$$\det(A_n-tB_n^T)=(-t)^n\det(1/tA_n-B_n^T)$$
we have $\det(A_n-tB_n^T)=(-t)^n(2t^{-n}-1)=(-1)^n(2-t^n)$.
Therefore, we have
$$\Delta_{K_{n,m}}(t)=(2t^n-1)(-1)^n(2-t^n)\doteq 2t^{n}-5+2t^{-n}.$$
Hence, since $d(K_{n,m})$ coincides with the genus of $\Sigma_{n,m}$,
we can see that the surface is the minimal Seifert surface.
Thus, we have $g(K_{n,m})=n$.

In the case of $K'_{n,m}$, we can do the similar computation to above by taking the corresponding generators in the Seifert surface.

The Seifert matrix $S'_n$ is $\begin{pmatrix}O_n&A'_n\\B'_n&O_n\end{pmatrix}$, where $O_n$ is the $n\times n$ zero matrix, $A'_n$ and $B'_n$ are $n\times n$ matrices satisfying the following:
$$A'_n=(a'_{ij}),\  a'_{ij}=\begin{cases}-2&j\ge i\\-3&j< i.\end{cases}\text{ and }B_n=(b_{ij}),\ b_{ij}=\begin{cases}-3&j\ge i\\-2&j< i.\end{cases}$$
Then we have
\begin{eqnarray*}
\Delta_{K'_{n,m}}&=&\det(tS'_n-{S'}_n^T)=\det\begin{pmatrix}O_n&tA'_n-{B'}_n^T\\tB'_n-{A'}_n^T&O_n\end{pmatrix}\\
&\doteq&6t^n-13+6t^{-n}.
\end{eqnarray*}

\hfill$\Box$
\end{proof}


\begin{lem}
Let $K_1$ and $K_2$ be two knots in two homology spheres $Y_1$, $Y_2$ respectively.
Then
$$g(K_1\#K_2)=g(K_1)+g(K_2)$$
holds.
\end{lem}
This is a classical result, however, we prove it here again.\\
\begin{proof}
Let $S\subset Y_1\#Y_2$ be the embedded separating sphere for $Y_1$ and $Y_2$.
We suppose that $S$ is separating $K_1\#K_2$ i.e., $(K_1\#K_2)\cap S$ are two points.
Let $\Sigma$ be the minimal genus Seifert surface of $K_1\#K_2$.
The set of the intersection $\Sigma\cap S$ consists of finite circles and single arc connecting the two points in the general position.
We take the inner most circle $C$ not including the arc in the interior.
$C$ bounds a disk in $\Sigma$ because $\Sigma$ is the minimal genus surface.
Then by cutting the disk and capping two new disks, we can decrease the number of the intersection circles.
We call the new embedded surface $\Sigma$ again.
This cut-and-past process preserves the genus of $\Sigma$.
The isotopy class of $\Sigma$ may be changed.
By iterating this process we vanish all the intersection circles.
Then $\Sigma=\Sigma_1\natural\Sigma_2$ is obtained and 
$$g(K_1\# K_2)=g(\Sigma)=g(\Sigma_1)+g(\Sigma_2)\ge g(K_1)+g(K_2)$$
holds.
Conversely, since $g(K_1)+g(K_2)\ge g(K_1\#K_2)$, we obtain
$g(K_1\#K_2)=g(K_1)+g(K_2)$.\hfill$\Box$
\end{proof}

Let $K$ be a knot in a homology sphere.
If $\Delta_K(t)$ cannot be decomposed as in $\Delta_K(t)=f_1(t)f_2(t)$ and $f_i(t)$ agrees with the Alexander polynomial of a knot in a homology sphere, then we call $\Delta_K(t)$ {\it A-irreducible}.
\begin{lem}
\label{prime}
Let $K$ be a knot in a homology sphere $Y$.
If $\Delta_K(t)$ is A-irreducible and $g(K)=d(K)$, then $K$ is prime. 
\end{lem}
\noindent{\begin{proof}
Suppose that $K$ is not prime.
Then $K$ is isotopic to a composite knot $K_1\#K_2$.
Then $\Delta_K(t)=\Delta_{K_1}(t)\Delta_{K_2}(t)$.
Hence, $d(K)=d(K_1)+d(K_2)$ holds.
Since $g(K)=d(K)$, we have $d(K)=g(K)=g(K_1)+g(K_2)\ge d(K_1)+d(K_2)$.
Therefore, $g(K_1)+g(K_2)= d(K_1)+d(K_2)$ holds.
From the inequalities $g(K_i)\ge d(K_i)$, $g(K_i)=d(K_i)$ holds for $i=1,2$.

On the other hand, since $K$ is A-irreducible, $\Delta_{K_1}(t)=1$ or $\Delta_{K_2}(t)=1$.
Thus $g(K_1)=0$ or $g(K_2)=0$ holds.
This means that $K$ is prime.
\hfill$\Box$
\end{proof}}
We prove $K_{n,m}$ is a prime knot.
\begin{lem}
$K_{n,m}$ and $K'_{n,m}$ are prime knots in $Y_{n,m}$, $Y'_{n,m}$ respectively.
\end{lem}
\noindent{\begin{proof}
The Alexander polynomials of $K_{n,m}$ and $K'_{n,m}$ are $2t^{n}-5+2t^{-n}$ and $6t^n-13+2t^{-n}$.
These polynomials are A-irreducible.
Because, the polynomials are completely decomposed as $\Delta_{K_{n,m}}\doteq 2t^{n}-5+2t^{-n}=(2t^n-1)(t^n-2)$ and $\Delta_{K'_{n,m}}\doteq 6t^{n}-13+6t^{-n}\doteq(2t^n-3)(3t^n-2)$ as a polynomial over ${\mathbb Z}$.
Any factor of these decompositions is an irreducible polynomial by Eisenstein's criterion
and is not an Alexander polynomial of a knot in a homology sphere.
Since the genus of $K_{n,m}$ and $K'_{n,m}$ is $n$, from Lemma~\ref{prime}, $K_{n,m}$ and $K'_{n,m}$ are prime.
\hfill$\Box$
\end{proof}}

\section{Boundary-sum irreducibility of $X({\bf x})_{n,m}$.}
Before proving Theorem~\ref{main}, we prove Lemma~\ref{equiv}.\\
\begin{proof}
Suppose that $X^4$ is boundary-sum reducible.
Then there exists a decomposition $X=X_1\natural X_2$ and $X_i$ is not homeomorphic to a 4-disk.
Suppose that either $\partial X_1$ or $\partial X_2$ is diffeomorphic to $S^3$.
We may assume $\partial X_1\cong S^3$.
Then $X$ is connected-sum $\hat{X}_1\#X_2$, where $\hat{X}_1$ is $X_1$ capped off by a 4-disk $D^4$ and $\hat{X}_1$ is not homeomorphic to $S^4$.
Then $X$ is irreducible.
Therefore we get the desired assertion.
\end{proof}
Here we prove Theorem~\ref{main}.\\
\begin{proof}
For any $\{\ast, 0\}$-sequence we set $X=X({\bf x})$.
The irreducible decomposition of $X_{n,m}$ is already done and unique.
Use Freedman's classification \cite{F} for the double of $X$.
Thus $X_{n,m}$ is irreducible.

We prove that $Y_{n,m}$ is a prime 3-manifold.
$Y_{n,m}$ is the $n$-fold cyclic branched cover of $Y(m)$ along $K_{n,m}$.
Namely, $Y_{n,m}/\langle \tau\rangle=Y(m)$, where $\tau=\tau_{n,m}$.
If $S\subset Y_{n,m}$ is an embedded 2-sphere, we assume that up to isotopy, $S$ satisfies either of the following conditions for any $g\in\langle \tau\rangle$ due to \cite{MSY} and \cite{Dun}:
\begin{itemize}
\item $g(S)\cap S=\emptyset$
\item $g(S)=S$.
\end{itemize}

Suppose that the first condition is satisfied.
$S$ does not intersect with the branch locus, namely, $S$ is projected to a sphere in $Y(m)$.
Since $Y(m)$ is a prime 3-manifold due to \cite{KOU}, then the sphere bounds a 3-ball in $Y(m)$.
Hence, lifting the ball to $Y_{n,m}$, we can find a 3-ball in $Y_{n,m}$ with the boundary $S$.

Suppose that the second condition is satisfied.
Then the action restricts on $S$.
The action is orientation-preserving, because if the action on $S$ is orientation-reversing,
then the quotient space has a connected-sum component of $L(2,1)$.
Then in the general position, $S$ transversely intersects with the branch locus at finite points.
By this argument, we can rule out the case where the branch locus is included in $S$.

This means that $\langle \tau\rangle$ acts on $S$ with the fixed points discrete.
The finite action of the 2-sphere is conjugate to the rotation in $SO(3)$ up to homotopy due to \cite{smale}.
In particular, the fixed points are two points.
Let $S'$ be an image of $S$ into $Y(m)$ and $S'\cap K_{n,m}$ are two points.
Since $Y(m)$ is prime, $S'$ bounds a 3-disk $D$ in $Y(m)$.
$D\cap K_{n,m}$ is the trivial arc, because $K_{n,m}$ is prime knot in $Y(m)$.
Since the branched cover along the trivial arc is a 3-disk, $S$ bounds a 3-disk in $Y_{n,m}$.

In each case, any embedded sphere in $Y_{n,m}$ bounds a 3-disk.
This means $Y_{n,m}$ is prime and it follows that $X_{n,m}$ is boundary-sum irreducible.

$Y'_{n,1}$ is the $n$-fold branched cover over $K_{n,1}$ in $Y'(1)=C(1)$.
This argument works for $Y'_{n,1}$.
This means $Y'_{n,1}$ is prime.
Thus, for any $n$, $(E_{n,1},\tau_{n,1})$ is boundary-sum irreducible ${\mathbb Z}_n$-order cork.
\hfill$\Box$
\end{proof}
Here, we put a quick proof of Theorem~\ref{irre}.\\
{\bf Proof of Theorem~\ref{irre}.}
It immediately follows from the latter of the proof.
\hfill$\Box$

\begin{figure}[htbp]
\begin{center}
\includegraphics{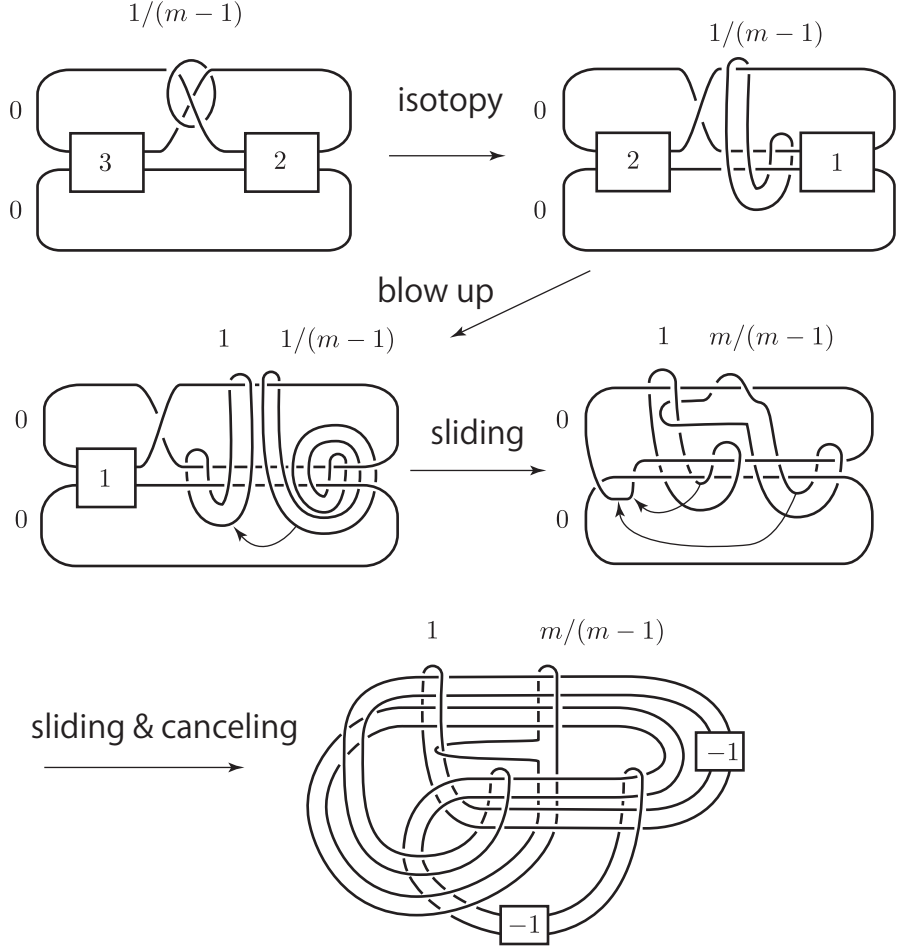}
\caption{The Dehn surgery diagram of $Y'(m)$.}
\label{YPnm}
\end{center}
\end{figure}

For any integer $m$ with $m\neq 1$, we do not know whether $E_{n,m}$ is boundary-sum irreducible or not.
We need to prove the primeness of $Y'(m)$.
The Dehn surgery diagram of $Y'(m)$ is drawn in {\sc Figure}~\ref{YPnm}.
This manifold is a Dehn surgery of $S^3_1(Pr(-3,3,-3))$.
\section{Proof of hyperbolicity.}
Finally, we prove Theorem~\ref{hyper}.\\
\begin{proof}
The output ``True" for the program HIKMOT means that the 3-manifold admits hyperbolic structure \cite{HIKMOT}.
To get True-output, we need apply Algorithm 2 in \cite{HIKMOT}.
The data after using the algorithm are updated in the site \cite{TM2}.
We can get ``True" for these four examples by running the data by HIKMOT.
\hfill$\Box$
\end{proof}

\end{document}